\documentclass[12pt,draft,leqno]{article}
\usepackage{amssymb, eucal, latexsym}
\newtheorem{theorem}{Theorem}[section]
\newtheorem{lm}[theorem]{Lemma}
\newtheorem{exa}[theorem]{Example}

\newtheorem{cor}[theorem]{Corollary}
\newtheorem{pro}[theorem]{Proposition}
\newtheorem{defi}[theorem]{Definition}

\newtheorem{nota}[theorem]{Notation}

\newtheorem{rem}[theorem]{Remark}

\newtheorem{fact}[theorem]{Fact}

\def\a{\alpha}
\def\b{\beta}

\def\UP{\Upsilon}

\def\sbe{\subseteq}

\def\nes{\neq\emptyset}

\def\ap{^\prime}

\def\st{\ |\ }

\def\AA{{\cal A}}
\def\BB{{\cal B}}

\def\MM{{\cal M}}
\def\OO{{\cal O}}
\def\PP{{\cal P}}
\def\TT{{\cal T}}

\def\UU{{\cal U}}

\def\ZZ{{\cal Z}}

\def\2{\mbox{{\bf 2}}}
\def\3{\mbox{{\bf 3}}}

\def\int{\mbox{{\rm int}}}

\def\doc{\hspace{-1cm}{\em Proof.}~~}
\def\sq{\hspace*{\fill} \hbox{\vrule\vbox{\hrule\phantom{o}\hrule}\vrule}}
\def\sqs{\sq \vspace{2mm}}

\def\RRRR{{\rm I}\!{\rm R}}
\def\NNNN{{\rm I}\!{\rm N}}

\title{{\LARGE\bf
Vietoris-type Topologies on Hyperspaces}
\vspace{0.35cm}
{\large\bf Elza Ivanova-Dimova}\\
\vspace{0.25cm}
 {\footnotesize Dept. of Math. and
Informatics, Sofia University,  J. Bourchier 5, 1126 Sofia,
Bulgaria}
}

\author{}

\date{}

\begin{document}
\maketitle
\begin{abstract}
{\footnotesize
\noindent We introduce a new
Vietoris-type hypertopology by means of the upper-Vietoris-type hypertopology  defined in \cite{DV} (it was called there {\em Tychonoff-type hypertopology}) and the lower-Vietoris-type hypertopology  introduced in \cite{ED}. We study this new Vietoris-type hypertopology and show that it is, in general, different from the Vietoris topology. Also, some of the results of E. Michael \cite{Mi} about hyperspaces with Vietoris topology are extended to analogous results for hyperspaces with Vietoris-type topology. We obtain as well some results about hyperspaces with Vietoris-type topology which concern some problems analogous to those regarded by H.-J. Schmidt in \cite{Sch}.}
\end{abstract}

{\footnotesize {\em  MSC:} 54B20, 54D10,  54C05.

}

\footnotetext[1]{{\footnotesize {\em E-mail address:}
elza@fmi.uni-sofia.bg}}

\baselineskip = \normalbaselineskip

\section*{Introduction}

In 1975, M. M. Choban \cite{Ch} introduced a new topology on the set of all closed
subsets of a topological space for obtaining a generalization of the famous Kolmogoroff
Theorem on operations on sets. This new topology is similar to the upper Vietoris
topology but is weaker than it. In 1998, G. Dimov and D. Vakarelov \cite{DV} defined another new hypertopology, which turned out to be a
generalized version of the Choban topology, and called it a {\em Tychonoff-type hypertopology} (it could be also called  an {\em upper-Vietoris-type hypertopology}\/); they used it for proving an isomorphism theorem for the
category of all Tarski consequence systems. Later on, the Tychonoff-type hypertopology was studied in details in \cite{DOT}. Further, in  \cite{ED}, a new lower-Vietoris-type hypertopology  was introduced and studied. In this talk we define a new
Vietoris-type hypertopology. We study it and show that, in general, it is different from the Vietoris topology. Also, some of the results of E. Michael \cite{Mi} about the hyperspaces with Vietoris topology are extended to analogous results for the hyperspaces with Vietoris-type topology. We obtain as well some results about the hyperspaces with Vietoris-type topology which concern some problems analogous to those regarded by H.-J. Schmidt in \cite{Sch}.

\section{Preliminaries}

We denote by $\NNNN$ the set of all natural numbers (hence, $0\not\in\NNNN$), by $\RRRR$ the real line (with its natural topology) and by $\overline{\RRRR}$ the set $\RRRR\cup\{-\infty,\infty\}$.

Let $X$ be a set. We denote by $|X|$ the cardinality of $X$ and by $\PP(X)$ (resp., by $\PP\ap(X)$) the set of all (non-empty) subsets of $X$. Let $\MM,\AA\subseteq\PP(X)$ and $A\subseteq X$. We set:

$\bullet$ $A_{\MM}^+:=\{M\in\MM\st M\subseteq A\}$;

$\bullet$ $\AA_{\MM}^+:=\{A_{\MM}^+\st A\in\AA\}$;

$\bullet$ $A^-_{\MM}:=\{M\in\MM\st M\cap A\not=\emptyset\}$

$\bullet$ $\AA_{\MM}^-:=\{A_{\MM}^-\st A\in\AA\}$;

$\bullet$ ${\it Fin}(X):=\{M\subseteq X\st 0<|M|<\aleph_0\}$;

$\bullet$ ${\it Fin_n}(X):=\{M\subseteq X\st 0<|M|\leq n\}$, where $n\in\NNNN$;

$\bullet$ $\displaystyle\AA^{\cap}:=\{\bigcap_{i=1}^k A_i\st k\in\NNNN, A_i\in\AA\}$.

$\bullet$ $\displaystyle\AA^{\cup}:=\{\bigcup_{i=1}^k A_i\st k\in\NNNN, A_i\in\AA\}$.

\noindent Let $(X,\TT)$ be a topological space. We put

$\bullet$ $CL(X):=\{M\subseteq X\st M \mbox{ is closed in }X,\ M\not=\emptyset\}$.

$\bullet$ ${\it Comp}(X):=\{M\in CL(X)\st M \mbox{ is compact}\}$.

\noindent When $\MM=CL(X)$, we will simply write $A^+$ and $A^-$ instead of $A_{\MM}^+$ and $A_{\MM}^-$; the same for subfamilies $\AA$ of $\PP(X)$. The closure of a subset $A$ of $X$ in $(X,\TT)$ will be denoted by $cl_{(X,\TT)}A$ or $\overline{A}^{(X,\TT)}$ (we will also write, for short, $cl_XA$ or $\overline{A}^X$ and even $\overline{A}^{\TT}$). By "neighborhood" we will mean an "open neighborhood". The regular spaces are not assumed to be $T_1$-spaces; by a $T_3$-space we mean a regular $T_1$-space.

If $X$ and $Y$ are sets and $f:X\longrightarrow Y$ is a function then, as usual, we denote by $f\upharpoonright X$ the function between  $X$ and $f(X)$ which is a restriction of $f$. If $(X,\TT)$ and $(Y,\OO)$ are topological spaces and $f:X\longrightarrow Y$ is an injection, then we say that $f$ is an {\em inversely continuous function} if the function $(f\upharpoonright X)^{-1}:f(X)\longrightarrow X$ is continuous.

Let $X$ be a topological space. Recall that the {\em  upper Vietoris topology} $\UP_{+X}$ on $CL(X)$ (called also {\em Tychonoff topology on} $CL(X)$) has as a base the family of all sets of the form $$U^+=\{F\in CL(X)\st F\subseteq U\},$$ where $U$ is open in $X$, and {\em the lower Vietoris topology} $\UP_{-X}$ on $CL(X)$ has as a subbase all sets of the form $$U^-=\{F\in CL(X)\st F\cap U\not=\emptyset\},$$ where $U$ is open in $X$. The {\em Vietoris topology} $\UP_X$ on $CL(X)$ is defined as the supremum of  $\UP_{+X}$ and $\UP_{-X}$, i.e., $\UP_{+X}\cup \UP_{-X}$ is a subbase for $\UP_X$.

\begin{defi}\label{Tichnoff type}{\rm (\cite{DV})}
\rm
(a) Let $(X,\TT)$ be a topological space and $\MM\subseteq\PP\ap(X)$. The topology $\UP_{+\MM}$ on $\MM$ having as a base the family $\TT_{\MM}^+$ is called a {\em Tychonoff topology on} $\MM$ generated by $(X,\TT)$. When $\MM=CL(X)$, then $\UP_{+\MM}$ is just the classical upper Vietoris topology $\UP_{+X}$ on $CL(X)$ (= Tychonoff topology $\UP_{+\MM}$ on $CL(X)$).

\medskip

(b) Let $X$ be a set and $\MM\subseteq\PP\ap(X)$. A topology $\OO$ on the set $\MM$ is called a {\em Tychonoff-type topology on} $\MM$ if the family $\OO\cap\PP(X)_{\MM}^+$ is a base for $\OO$.
\end{defi}

A Tychonoff topology on $\MM$ is always a Tychonoff-type topology on $\MM$, but not viceversa (see \cite{DOT}).

\begin{fact}\label{u-Vietoris-type}{\rm (\cite{DV})}
\rm
Let $X$ be a set, $\MM\subseteq\PP\ap(X)$ and $\OO$ be a topology on $\MM$. Then the family $$\BB_{\OO}:=\{A\subseteq X\st A_{\MM}^+\in\OO\}$$ contains $X$ and is closed under
finite intersections; hence, it can serve as a base for a topology $$\TT_{+\OO}$$ on $X$. When $\OO$ is a Tychonoff-type topology, the family $(\BB_{\OO})_{\MM}^+$ is a base for $\OO$.
\end{fact}

\begin{defi}\label{2.5new}{\rm (\cite{DV})}
\rm
Let $X$ be a set, $\MM\subseteq\PP\ap(X)$ and $\OO$ be a topology on $\MM$. Then we  say that the topology $\TT_{+\OO}$ on $X$, introduced in Fact \ref{u-Vietoris-type}, is {\em the plus-topology on $X$  induced by the topological space} $(\MM,\OO)$.
\end{defi}

\begin{defi}\label{LVT}{\rm (\cite{ED})}
\rm
Let $(X,\TT)$ be a topological space and $\MM\subseteq\PP\ap(X)$. The topology $\UP_{-\MM}$ on $\MM$ having as a subbase the family $\TT_{\MM}^-$ is called a {\em lower Vietoris topology on} $\MM$ generated by $(X,\TT)$. When $\MM=CL(X)$, then $\UP_{-\MM}$ is just the classical lower Vietoris topology $\UP_{-X}$ on $CL(X)$.
\end{defi}

\begin{defi}\label{lower-Vietoris-type}{\rm (\cite{ED})}
\rm
Let $X$ be a set, $\MM\subseteq\PP\ap(X)$, $\OO$ be a topology on $\MM$. We say that $\OO$ is a {\em lower-Vietoris-type topology on} $\MM$, if $\OO\cap\{A^-_{\MM}\st A\subseteq X\}$ is a subbase for $\OO$.
\end{defi}

A lower Vietoris topology on $\MM$ is always a lower-Vietoris-type topology on $\MM$, but not viceversa (see \cite{ED}).

\begin{fact}\label{Vietoris-type}{\rm (\cite{ED})}
\rm
Let $X$ be a set, $\MM\subseteq\PP\ap(X)$ and $\OO$ be a topology on $\MM$. Then the family $$\PP_{\OO}:=\{A\subseteq X\st A_{\MM}^-\in\OO\}$$ contains $X$, and, hence, can serve as a subbase for a topology $$\TT_{-\OO}$$ on $X$. The family $\PP_{\OO}$ is closed under arbitrary unions. When $\OO$ is a lower-Vietoris-type topology, the family $(\PP_{\OO})_{\MM}^-$ is a subbase for $\OO$.
\end{fact}

\begin{defi}\label{2.5}
\rm
Let $X$ be a set, $\MM\subseteq\PP\ap(X)$ and $\OO$ be a topology on $\MM$. Then we  say that the topology $\TT_{-\OO}$ on $X$, introduced in Fact \ref{Vietoris-type}, is {\em the minus-topology on $X$ induced by the topological space} $(\MM,\OO)$.
\end{defi}

All undefined here notions and notation can be found in \cite{AP,E2}.

\section{The Vietoris-type topologies on hyperspaces}

\begin{defi}\label{Vie}
\rm
Let $(X,\TT)$ be a topological space and $\MM\subseteq\PP\ap(X)$. The topology $\UP_{\MM}$ on $\MM$ having as a subbase the family $\TT_{\MM}^+\cup\TT_{\MM}^-$ is called {\em the Vietoris topology on} $\MM$. When $\MM=CL(X)$ then $\UP_{\MM}\equiv\UP_X$.
\end{defi}

\begin{defi}\label{vietoris-type}
\rm
Let $X$ be a set, $\MM\subseteq\PP\ap(X)$ and $\OO$ be a topology on $\MM$. The topology $\OO$ is called a {\em Vietoris-type topology on} $\MM$ if the family $(\BB_{\OO})_{\MM}^+\cup(\PP_{\OO})_{\MM}^-$ is a subbase for $\OO$.
\end{defi}

\begin{rem}\label{2.1'}
\rm
Let $(X,\TT)$ be a topological space and $\MM\subseteq\PP\ap(X)$. Then, clearly, $\UP_{\MM}$ is a Vietoris-type topology on $\MM$. As we will see in Example \ref{novietoris} below, the converse is not always true even when $\MM=CL(X)$. Note that when $\MM\subseteq CL(X)$, then $\UP_{\MM}\equiv(\UP_X)|_{\MM}$.
\end{rem}

\begin{nota}\label{V-top}
\rm
Let $X$ be a set, $\MM\subseteq\PP\ap(X)$ and $\OO$ be a topology on $\MM$. We denote by $$\TT_{\OO}$$ the topology on $X$ having $\BB_{\OO}\cup\PP_{\OO}$ as a subbase and we say that $\TT_{\OO}$ is the {\em V-topology on} $X$ {\em induced by the topological space} $(\MM,\OO)$. We denote by $$\OO_u$$ the topology on $\MM$ having $(\BB_{\OO})_{\MM}^+$ as a base, and by $$\OO_l$$ the topology on $\MM$ having $(\PP_{\OO})_{\MM}^-$ as a subbase.
\end{nota}

\begin{defi}\label{natural}
\rm
Let $X$ be a set and $\MM\subseteq\PP(X)$. We say that $\MM$ is a {\em natural family in} $X$ if $\{x\}\in\MM$ for any $x\in X$.
\end{defi}

The following assertion can be easily proved.

\begin{fact}\label{supremum}
Let $X$ be a set, $\MM\subseteq\PP\ap(X)$ and $\OO$ be a topology on $\MM$. Then $\TT_{\OO}=\TT_{+\OO}\vee\TT_{-\OO}$, $\OO_u$ is a Tychonoff-type topology on $\MM$ and $\OO_l$ is a lower-Vietoris-type topology on $\MM$. $\OO$ is a Vietoris-type topology on $\MM$ iff $\OO=\OO_u\vee\OO_l$.
\end{fact}

\begin{fact}\label{fv3}
Let $(X,\TT)$ be a topological space. If $\MM\subseteq CL(X)$, then $(\UP_{+X})|_{\MM}$ and $(\UP_{-X})|_{\MM}$ are Vietoris-type topologies on $\MM$. Moreover, if $\MM\subseteq\PP\ap(X)$, then $\UP_{+\MM}$ and $\UP_{-\MM}$ are Vietoris-type topologies on $\MM$.
\end{fact}

\doc Let $\MM\subseteq\PP\ap(X)$. Set $\OO=\UP_{+\MM}$. Then $\TT_{\MM}^+$ is a base for $\OO$. Since $\TT_{\MM}^+\subseteq(\BB_{\OO})_{\MM}^+\subseteq(\BB_{\OO})_{\MM}^+\cup(\PP_{\OO})_{\MM}^-\subseteq\OO$, we get that $(\BB_{\OO})_{\MM}^+\cup(\PP_{\OO})_{\MM}^-$ is a (sub)base for $\OO$. Hence $\UP_{+\MM}$ is a Vietoris-type topology on $\MM$. Analogously, we get that $\UP_{-\MM}$ is a Vietoris-type topology on $\MM$.

Let now $\MM\subseteq CL(X,\TT)$. $\TT^+$ is a base for $\UP_{+X}$. Hence $(\TT^+)|_{\MM}=\{U^+\cap\MM\st U\in\TT\}$ is a base for $(\UP_{+X})|_{\MM}$. Since $U^+\cap\MM=U_{\MM}^+$, we get that $\TT_{\MM}^+$ is a base for $(\UP_{+X})|_{\MM}$. Hence $\UP_{+\MM}\equiv(\UP_{+X})|_{\MM}$ and thus $(\UP_{+X})|_{\MM}$ is a Vietoris-type topology on $\MM$. Analogously, we get that $(\UP_{-X})|_{\MM}$ is a Vietoris-type topology on $\MM$ (and $(\UP_{-X})|_{\MM}\equiv\UP_{-\MM}$). \sqs

\begin{defi}\label{strongVietoris}
\rm
Let $X$ be a set, $\MM\subseteq\PP\ap(X)$ and $\OO$ be a Vietoris-type topology on $\MM$. Then $\OO$ is called a {\em strong Vietoris-type topology on} $\MM$ if $\TT_{+\OO}\equiv\TT_{-\OO}$.
\end{defi}

\begin{pro}\label{strV}
Let $(X,\TT)$ be a topological space, $\MM\subseteq CL(X)$ and $\MM$ be a natural family. Let $\OO=(\UP_X)|_{\MM}$. Then $\OO$ is a strong Vietoris-type topology on $\MM$ and $\TT\equiv\TT_{\OO}$. In particular, for every $T_1$-space $X$, $\UP_X$ is a strong Vietoris-type topology on $CL(X)$.
\end{pro}

\doc By Remark \ref{2.1'}, $\OO=(\UP_X)|_{\MM}$ is a Vietoris-type topology on $\MM$. Obviously, $\TT\subseteq\BB_{\OO}\cap\PP_{\OO}$. We will show that $\BB_{\OO}\cup\PP_{\OO}\subseteq\TT$. Clearly, this will imply that $\TT_{+\OO}=\TT_{-\OO}=\TT=\TT_{\OO}$. Let $A\subseteq X$ and $A_{\MM}^+\in\OO$. It is well-known that the family $\{<U_1,\dots,U_n>\st n\in\NNNN,\ U_i\in\TT\mbox{ for } i\in\{1,\dots,n\}\}$, where $\displaystyle<U_1,\dots,U_n>=(\bigcup_{i=1}^nU_i)^+\cap\bigcap_{i=1}^n(U_i)^-$, is a base for $\UP_X$. Let $x\in A$. Since $\MM$ is a natural family, we get that $\{x\}\in A_{\MM}^+$. Thus there exist $U_1,\dots,U_n\in\TT$ such that $\displaystyle\{x\}\in(\bigcup_{i=1}^nU_i)_{\MM}^+\cap\bigcap_{i=1}^n(U_i)_{\MM}^-\subseteq A_{\MM}^+$. Then $\displaystyle x\in U=\bigcap_{i=1}^nU_i$. Using again naturality of $\MM$, we obtain that $U\subseteq A$. Hence, $A\in\TT$. So, $\BB_{\OO}\subseteq\TT$. Let now $A\subseteq X$ and $A_{\MM}^-\in\OO$. Then, arguing as above, we get that $A\in\TT$. Thus, $\PP_{\OO}\subseteq\TT$. Hence $\BB_{\OO}\cup\PP_{\OO}\subseteq\TT$. \sqs

\begin{exa}\label{noVt}
\rm
There exists a $T_0$-space $(X,\TT)$ such that $\OO=\UP_X$ is not a strong Vietoris-type topology on $CL(X)$. Also, we have that $\TT\not=\TT_{-\OO}$, $\TT\not=\TT_{+\OO}$, $\TT\not=\TT_{\OO}$ and $\TT_{\OO}=\TT_{-\OO}$.
\end{exa}

\doc Let $X=\{0,1,2\}$ and $\TT=\{\emptyset,X,\{0\},\{0,2\}\}$. Then $(X,\TT)$ is a $T_0$-space and $$CL(X)=\{X,\{1\},\{1,2\}\}.$$ The topology $\UP_{-X}$ has as a subbase the family $$\{\emptyset^-, X^-,\{0\}^-,\{0,2\}^-\}=\{\emptyset,CL(X),\{X\},\{\{1,2\},X\}\}.$$ Thus $\UP_{-X}=\{\emptyset,CL(X),\{X\},\{\{1,2\},X\}\}$. The topology $\UP_{+X}$ has as a base the family $\{\emptyset^+,X^+,\{0\}^+,\{0,2\}^+\}=\{\emptyset,CL(X)\}$. Thus $\UP_{+X}=\{\emptyset,CL(X)\}$. Hence $\UP_X=\UP_{+X}\vee\UP_{-X}=\UP_{-X}$. Set $\OO=\UP_{X}$. Then $\BB_{\OO}=\{A\subseteq X\st A^+\in\OO\}$ and $\PP_{\OO}=\{A\subseteq X\st A^-\in\OO\}$. We have that $\emptyset^+=\emptyset^-=\emptyset$, $X^+=X^-=CL(X)$, $\{0\}^-=\{X\}$, $\{0\}^+=\emptyset$, $\{1\}^+=\{\{1\}\}$, $\{1\}^-=CL(X)$, $\{2\}^+=\emptyset$, $\{2\}^-=\{X,\{1,2\}\}$, $\{0,1\}^+=\{\{1\}\}$, $\{0,1\}^-=CL(X)$, $\{0,2\}^+=\emptyset$, $\{0,2\}^-=\{X,\{1,2\}\}$, $\{1,2\}^+=\{\{1\},\{1,2\}\}$, $\{1,2\}^-=CL(X)$. Hence $\BB_{\OO}=\{\emptyset,X,\{0\},\{2\},\{0,2\}\}$ and $\PP_{\OO}=\{\emptyset,X,\{0\},\{1\},\{2\},\{0,1\},\{0,2\},\{1,2\}\}$. Thus $\TT_{+\OO}=\BB_{\OO}$ and $\TT_{-\OO}=\PP_{\OO}$. Therefore $\TT_{+\OO}\not=\TT_{-\OO}$. Hence $\UP_X$ is not a strong Vietoris-type topology on $CL(X)$. Note that $\TT\not=\TT_{+\OO}$, $\TT\not=\TT_{-\OO}$ and $\TT\not=\TT_{\OO}$. Since  $\BB_{\OO}\subseteq\PP_{\OO}$, we get that $\TT_{+\OO}\subseteq\TT_{-\OO}$ and thus $\TT_{\OO}=\TT_{-\OO}$. \sqs

\begin{exa}\label{noVt1}
\rm
There exists a $T_0$-space $(X,\TT)$ such that $\OO=\UP_X$ is not a strong Vietoris-type topology on $CL(X)$ and $\TT_{\OO}=\TT_{-\OO}\supsetneqq\TT_{+\OO}=\TT$.
\end{exa}

\doc Let $X=\{0,1,2\}$ and $\TT=\{\emptyset,X,\{0\},\{0,2\},\{1,2\},\{2\}\}$. Then $(X,\TT)$ is a $T_0$-space and $CL(X)=\{X,\{1,2\},\{1\},\{0\},\{0,1\}\}$. The topology $\UP_{-X}$ has as a subbase the family
$$\begin{array}{l}
\{\emptyset^-,X^-,\{0\}^-,\{0,2\}^-,\{1,2\}^-,\{2\}^-\}=\\ \{\emptyset,CL(X),\{X,\{0\},\{0,1\}\},\{X,\{0\}, \{0,1\},\{1,2\}\},\{\{1\},\{0,1\},X,\{1,2\}\},\\ \{X,\{1,2\}\}\}.
\end{array}$$
Then $\UP_{-X}=\{\emptyset,CL(X),\{X,\{0\},\{0,1\}\},\{X,\{0\},\{0,1\},\{1,2\}\},\\ \{X,\{1\},\{0,1\},\{1,2\}\},\{X,\{1,2\}\}, \{X,\{0,1\}\},\{X\},\{X,\{0,1\},\{1,2\}\}\}$. The topo\-logy $\UP_{+X}$ has as a base the family $$\{\emptyset^+,X^+,\{0\}^+,\{0,2\}^+,\{1,2\}^+,\{2\}^+\}=\{\emptyset,CL(X),\{\{0\}\},\{\{1\},\{1,2\}\}\}.$$ Hence $\UP_{+X}=\{\emptyset,CL(X),\{\{0\}\},\{\{1\},\{1,2\}\},\{\{0\},\{1\},\{1,2\}\}\}$ and $\OO=\UP_X=\UP_{-X}\vee\UP_{+X}=\UP_{-X}\cup\UP_{+X}\cup\{\{1,2\}\}\cup\{\{0\},\{1,2\}\}\cup\{X,\{0\},\{1,2\}\}\cup \{X,\{0\}\}\cup\{X,\{1\},\{1,2\}\}\cup\{X,\{0\},\{1\},\{1,2\}\}$. We have that $\emptyset^+=\emptyset^-=\emptyset$, $X^+=X^-=CL(X)$, $\{0\}^-=\{X,\{0\},\{0,1\}\}$, $\{0\}^+=\{\{0\}\}$, $\{1\}^-=\{X,\{1\},\{0,1\},\{1,2\}\}$, $\{1\}^+=\{\{1\}\}$, $\{2\}^-=\{X,\{1,2\}\}$, $\{2\}^+=\emptyset$, $\{0,1\}^-=CL(X)$, $\{0,1\}^+=\{\{0\},\{1\},\{0,1\}\}$, $\{0,2\}^-=\{X,\{0\},\{0,1\},\{1,2\}\}$, $\{0,2\}^+=\{\{0\}\}$, $\{1,2\}^-=\{X,\{1\},\{0,1\},\{1,2\}\}$, $\{1,2\}^+=\{\{1\},\{1,2\}\}$. Hence $$\BB_{\OO}=\{\emptyset,X,\{0\},\{2\},\{0,2\},\{1,2\}\}=\TT_{+\OO}=\TT$$ and $\PP_{\OO}=\{\emptyset,X,\{0\},\{2\},\{1,2\},\{0,2\},\{1\},\{0,1\}\}=\TT_{-\OO}\not=\TT$. Thus $$\TT_{\OO}=\TT_{-\OO}\supsetneqq\TT_{+\OO}\equiv\TT.$$ \sqs

\begin{exa}\label{novietoris}
\rm
Let $\UU=\{(\alpha,\beta)\st\alpha,\beta\in\overline{\RRRR}\}$ and $\RRRR$ be regarded with its natural topology. Then
the topology $\OO$ on $CL(\RRRR)$ having as a subbase the family $\UU^-\cup\UU^+$ is a strong Vietoris-type topology on $CL(\RRRR)$  different from the Vietoris topology
$\UP_{\RRRR}$ on $CL(\RRRR)$. Also, $\TT_{\OO}$ coincides with the natural topology $\TT$ on $\RRRR$.
\end{exa}

\doc Clearly, $\UU\subseteq\BB_{\OO}\cap\PP_{\OO}$ and thus $(\BB_{\OO})^+\cup(\PP_{\OO})^-$ is a subbase for $\OO$. Hence $\OO$ is a Vietoris-type topology on $CL(\RRRR)$.
We will show that $\TT_{+\OO}\equiv\TT_{-\OO}\equiv\TT$. Then, in particular, we will obtain that $\OO$ is a strong Vietoris-type topology on $CL(\RRRR)$. Note first that $\UU^\cap=\UU$ and $(\UU^+)^\cap=\UU^+$.
Fufther, it is not difficult to see (using the fact that $CL(\RRRR)$ is a natural family) that $\PP_{\OO}\cup\BB_{\OO}\sbe\TT$ (see Fact \ref{Vietoris-type} and Fact \ref{u-Vietoris-type} for the notation $\PP_\OO$ and $\BB_\OO$). Since $\UU$ is a base for $\TT$ and $\UU\subseteq\BB_{\OO}\cap\PP_{\OO}$, we  get that $\TT_{+\OO}\equiv\TT_{-\OO}\equiv\TT$. So, $\OO$ is a strong Vietoris-type topology on $CL(\RRRR)$ and $\TT_\OO=\TT$.

For showing that $\OO\neq\UP_{\RRRR}$, we will prove that $((0,1)\cup(1,2))^+\not\in\OO$. Let $F=\{\frac{1}{2},\frac{3}{2}\}$. Then $F\in((0,1)\cup(1,2))^+$.  We will show that if $n\in\NNNN$, $U_1,\dots,U_n, V\in\UU$
and $\displaystyle F\in V^+\cap\bigcap_{i=1}^nU_i^-=O$, then $O\not\subseteq((0,1)\cup(1,2))^+$. Indeed, there exist $\a, \b\in\overline{\RRRR}$ such that $\a<\b$ and $V=(\a,\b)$.
Then $F\subset(\alpha,\beta)$ and $F\cap U_i\not=\emptyset$, for $i=1,\ldots, n$.
Clearly, $\alpha<\frac{1}{2}<1<\frac{3}{2}<\beta$. Let $G=[\frac{1}{2},\frac{3}{2}]$. Then $G\in CL(\RRRR)$, $F\subset G$, $G\subset(\alpha,\beta)$ and $G\cap U_i\nes$, for all $i=1,\ldots,n$. Therefore $G\in O$, but $G\not\in((0,1)\cup(1,2))^+$ since $1\in G$ and $1\not\in(0,1)\cup(1,2)$.
So, $((0,1)\cup(1,2))^+\in \UP_{\RRRR}$ but $((0,1)\cup(1,2))^+\not\in\OO$. Hence $\UP_{\RRRR}\not=\OO$ (and, clearly, $\OO\subset \UP_{\RRRR}$). \sqs

\section{Some properties of the hyperspaces with Vieto\-ris-type topologies}

In this section, some of the results of E. Michael \cite{Mi} concerning hyperspaces with Vietoris topology will be extended to analogous results for the hyperspaces with Vietoris-type topology.

\begin{pro}\label{2.7.20a}
Let $(X,\TT)$ be a space, $\MM\subseteq\PP\ap(X)$, $n\in\NNNN$, $Fin_n(X)\subseteq\MM$, $\OO$ be a Vietoris-type topology on $\MM$ and $\TT_{\OO}\subseteq\TT$. Let $J_n(X)$ be the subspace of $(\MM,\OO)$ consisting of all sets of cardinality $\leq n$. Then the map $j_n:X^n\longrightarrow J_n(X)$, where $j_n(x_1,\dots,x_n)=\{x_1,\dots,x_n\}$, is continuous.
\end{pro}

\doc Let $x=(x_1,\dots,x_n)\in X^n$ and $\displaystyle j_n(x)\in U_{\MM}^+\cap\bigcap_{i=1}^k(U_i)_{\MM}^-=O$, where $U\in\BB_{\OO}$, $k\in\NNNN$ and $U_i\in\PP_{\OO}$, $\forall i=1,\dots,k$. Then $\{x_1,\dots,x_n\}\subseteq U$ and $\{x_1,\dots,x_n\}\cap U_i\not=\emptyset$, $\forall i=1,\dots,k$. Hence, $\forall i\in\{1,\dots,k\}$ $\exists s(i)\in\{1,\dots,n\}$ such that $x_{s(i)}\in U_i$. Set $\displaystyle V=U^n\cap\prod_{t=1}^nV_t$, where
$V_t=\left\{\begin{array}{l}U_i,\mbox{ if }t=s(i)\mbox{ for some }i\in\{1,\dots,k\}\\
X,\mbox{ otherwise }\end{array}\right.$. Then, clearly, $x\in V$ and $j_n(V)\subseteq O$. So, $j_n$ is a continuous function. \sqs

\begin{pro}\label{2.7.20b}
Let $(X,\TT)$ be a space, $\MM\subseteq\PP\ap(X)$ be a natural family, $\OO$ be a Vietoris-type topology on $\MM$ and $\TT_{\OO}=\TT$. Then $j_1:(X,\TT)\longrightarrow J_1(X)$ is a homeomorphism.
\end{pro}

\doc Using Proposition \ref{2.7.20a}, we get that $j_1$ is a continuous bijection. We will prove that $j_1^{-1}$ is continuous. Let $\{x\}\in J_1(X)$, $U\in\BB_{\OO}$, $k\in\NNNN$, $U_1,\dots,U_k\in\PP_{\OO}$ and $\displaystyle x\in V=U\cap\bigcap_{i=1}^kU_i$. Then $\displaystyle\{x\}\in O=U_{\MM}^+\cap\bigcap_{i=1}^k(U_i)_{\MM}^-$ and, clearly, $(j_1)^{-1}(O)\subseteq V$. \sqs

\begin{pro}\label{2.7.20v}
If $(X,\TT)$ is a $T_2$-space, $\MM\subseteq\PP'(X)$, $\MM$ is a natural family, $\OO$ is a Vietoris-type topology on $\MM$ and $\TT_{-\OO}\supseteq\TT$, then $J_1(X)$ is closed in $(\MM,\OO)$.
\end{pro}

\doc Let $M\in\MM\setminus J_1(X)$. Then there exist $x,y\in M$ such that $x\not=y$. Since $(X,\TT)$ is a Hausdorff space, there exist $U_1,\dots,U_k, V_1,\dots,V_l\in\PP_{\OO}$ (where $k,l\in\NNNN$) such that $\displaystyle x\in U=\bigcap_{i=1}^kU_i$, $\displaystyle y\in V=\bigcap_{j=1}^lV_j$ and $U\cap V=\emptyset$. Set $\displaystyle O=\bigcap_{i=1}^k(U_i)_{\MM}^-\cap\bigcap_{j=1}^l(V_j)_{\MM}^-$. Then $M\in O$ and $O\cap J_1(X)=\emptyset$. Hence $J_1(X)$ is a closed subset of $(\MM,\OO)$. \sqs

\begin{pro}\label{2.7.20g}
Let $(X,\TT)$ be a space, $\MM\subseteq\PP'(X)$, $Fin_2(X)\subseteq\MM$, $\OO$ be a Vietoris-type topology on $\MM$, $\TT_{\OO}\subseteq\TT$ and $J_1(X)$ be closed in $(\MM,\OO)$. Then $X$ is a $T_2$-space.
\end{pro}

\doc Let $x_1,x_2\in X$, $x_1\not=x_2$. Then $\{x_1,x_2\}\in\MM$ and $\{x_1,x_2\}\not\in J_1(X)$. Since $J_1(X)$ is closed, $\exists V\in\BB_{\OO}$ and $U_1,\dots,U_k\in\PP_{\OO}$ such that $\displaystyle\{x_1,x_2\}\in V_{\MM}^+\cap\bigcap_{i=1}^n(U_i)_{\MM}^-=O$ and $O\cap J_1(X)=\emptyset$. Then $\{x_1,x_2\}\subset V$ and $\{x_1,x_2\}\cap U_i\not=\emptyset$ $\forall i\in\{1,\dots,k\}$. Let $\UU=\{U_1,\dots,U_k\}$, $V_1=\bigcap\{U\in\UU\st x_1\in U\}$ and $V_2=\bigcap\{U\in\UU\st x_2\in U\}$. Then $V_1,V_2\in\TT$. We will show that $x_1\in V_1$, $x_2\in V_2$ and $V\cap V_1\cap V_2=\emptyset$. Indeed, since $\{x_1\}\not\in O$, $\exists U\in\UU$ such that $x_1\not\in U$. Then $x_2\in U$ (because $U\cap\{x_1,x_2\}\not=\emptyset$). Hence $V_2\not=\emptyset$ and, obviously, $x_2\in V_2$. Analogously, we get that $x_1\in V_1$. Suppose that $\exists x\in V\cap V_1\cap V_2$. Then $x\in U_i$ for every $i=1,dots,k$. (Indeed, for every $U\in\UU$, we have that either $x_1\in U$ or $x_2\in U$.) Since $x\in V$, we get that $\{x\}\in O\cap J_1(X)$, a contradiction. Hence $x_1\in V\cap V_1$, $x_2\in V\cap V_2$ and $(V\cap V_1)\cap(V\cap V_2)=\emptyset$. \sqs

\begin{pro}\label{2.7.20d}
Let $X$ be a set, $\MM\subseteq\PP'(X)$, $Fin(X)\subseteq\MM$ and $\OO$ be a Vietoris-type topology on $\MM$. Then $\displaystyle J(X)=\bigcup_{i\in\NNNN}J_i(X)$ is dense in $(\MM,\OO)$.
\end{pro}

\doc Let $n\in\NNNN$, $U\in\BB_{\OO}\setminus\{\emptyset\}$, $U_i\in\PP_{\OO}\setminus\{\emptyset\}$ for $i=1,\dots,n$ and $\displaystyle O=U_{\MM}^+\cap\bigcap_{i=1}^n(U_i)_{\MM}^-\not=\emptyset$. Then $O$ is a basic open subset of $(\MM,\OO)$ and it is enough to show that $O\cap J(X)\not=\emptyset$. Since $O\not=\emptyset$, we have that $U\cap U_i\not=\emptyset$ for every $i\in\{1,\dots,n\}$. Indeed, if $M\in O$ then $M\subseteq U$ and $M\cap U_i\not=\emptyset$ $\forall i=1,\dots,n$. So, $\forall i=1,\dots,n$, $\exists x_i\in U\cap U_i$. Then $\{x_1,\dots,x_n\}\in J(X)\cap O$. Hence, $J(X)$ is dense in $(\MM,\OO)$. \sqs

\begin{defi}\label{predteglo}
\rm
Let $X$ be a set and $\PP\subseteq\PP(X)$. Set $w(X,\PP)=\min\{|\PP'|\st(\PP'\subseteq\PP)\wedge(\forall U\in\PP\mbox{ and }\forall x\in U\ \exists V\in\PP'\mbox{ such that }x\in V\subseteq U)\}$.
\end{defi}

Clearly, $w(X,\PP)\leq|\PP|$; also, when $\PP$ is a topology on $X$, then $w(X,\PP)$ is just the weight of the topological space $(X,\PP)$.

\begin{fact}\label{2.21}
Let $(X,\TT)$ be a topological space and $\PP$ be a subbase for $(X,\TT)$. Then:

\smallskip

\noindent (a) the families $\PP\ap$ from Definition \ref{predteglo} are also subbases for $(X,\TT)$;

\smallskip

\noindent (b) if $w(X,\TT)\geq\aleph_0$ then $w(X,\PP)\geq w(X,\TT)$.
\end{fact}

\begin{rem}\label{2.21'}
\rm
Let $X$ be a set and $\PP\subseteq\PP(X)$. Then, clearly, there exists a unique topology $\TT(\PP)$ on $X$ for which $\PP\cup\{X\}$ is a subbase. Obviously, if $\bigcup\PP=X$ then $\PP$ is a subbase for $\TT(\PP)$. Hence, in Definition \ref{predteglo} we can always assume that $X$ is a topological space and $\PP\cup\{X\}$ is a subbase for $X$.
\end{rem}

In connection with Fact \ref{2.21}, note that the following assertion holds (it should be well-known):

\begin{lm}\label{teglo}
Let $X$ be a space, $w(X)=\tau\geq\aleph_0$ and $\PP$ be a subbase for $X$. Then  there exists a $\PP'\subseteq\PP$, such that $|\PP'|=\tau$ and $\PP\ap$ is a subbase for $X$.
\end{lm}

\doc Let $\BB=\PP^{\cap}$. Then there exists a base $\BB'\subseteq\BB$ for $X$ such that $|\BB'|=\tau$. For every element $U$ of $\BB\ap$ fix a finite subfamily $\UU(U)$ of $\PP$ such that $U=\bigcap\UU(U)$. Set $\displaystyle\PP'=\{U'\in\PP\st\exists U\in\BB'\mbox{ such that }U'\in\UU(U)\}$. Then $\PP'$ is a subbase for $X$ (because $\BB'\subseteq(\PP')^{\cap}$ and $\BB\ap$ is a base for $X$) and $|\PP'|\leq|\BB'|\cdot\aleph_0=|\BB'|=\tau$. Since $|(\PP')^{\cap}|=|\PP'|$ and $(\PP')^{\cap}$ is a base for $X$, we get that $|\PP\ap|\geq\tau$. Hence $|\PP'|=\tau$. \sqs

\begin{pro}\label{3.3'}
Let $X$ be a set, $\MM\subseteq\PP\ap(X)$ and $\OO$ be a lower-Vietoris-type topology on $\MM$. Let $\tau\geq\aleph_0$, $\PP'\subseteq\PP(X)$ and $(\PP')_{\MM}^-$ be a subbase for $\OO$. If $w(X,\PP')\leq\tau$, then $w(\MM,\OO)\leq\tau$. In particular, if $\tau\geq\aleph_0$ and $w(X,\PP_{\OO})\leq\tau$ then $w(\MM,\OO)\leq\tau$.
\end{pro}

\doc If $w(\MM,\OO)<\aleph_0$ then $w(\MM,\OO)\leq\tau$. Hence, we can suppose that $w(\MM,\OO)\geq\aleph_0$. Let $\PP\subseteq\PP\ap$, $|\PP|\leq\tau$ and for every $U\in\PP\ap$ and for every $x\in U$ there exists a $V\in\PP$, such that $x\in V\subseteq U$. We will prove that $(\PP_{\MM}^-)^{\cap}$ is a base for $(\MM,\OO)$. Indeed, let $M\in\MM$, $O\in\OO$ and $M\in O$. Then there exist $U_1,\dots,U_n\in\PP\ap$ such that $\displaystyle M\in\bigcap_{i=1}^n(U_i)_{\MM}^-\subset O$. Hence $M\cap U_i\not=\emptyset$ for every $i=1,\dots,n$. Let $x_i\in U_i\cap M$. Then, for every $i=1,\dots,n$, there exists $V_i\in\PP$ such that $x_i\in V_i\subseteq U_i$. Hence $\displaystyle M\in\bigcap_{i=1}^n(V_i)_{\MM}^-\subset\bigcap_{i=1}^n(U_i)_{\MM}^-\subset O$. So, $(\PP_{\MM}^-)^{\cap}$ is a base for $(\MM,\OO)$ and thus $|(\PP_{\MM}^-)^{\cap}|\geq\aleph_0$. Therefore $|\PP_{\MM}^-|=|(\PP_{\MM}^-)^{\cap}|\geq\aleph_0$ and we obtain that $w(\MM,\OO)\leq|(\PP_{\MM}^-)^{\cap}|\leq|\PP|\leq\tau$. \sqs

\begin{pro}\label{teglovi}
Let $X$ be a set, $\MM\subseteq\PP\ap(X)$, $\OO$ be a Vietoris-type topology on $\MM$. Let $\tau\geq\aleph_0$, $w(X,\PP_{\OO})\leq\tau$ and there exists a family $\BB\subseteq\BB_{\OO}$ such that $|\BB|\leq\tau$ and $\forall M\in\MM$, $\forall U\in\BB_{\OO}$ with $M\subseteq U$, $\exists V\in\BB$ with $M\subseteq V\subseteq U$. Then $w(\MM,\OO)\leq\tau$.
\end{pro}

\doc Clearly, Proposition \ref{3.3'} implies that $w(\MM,\OO_l)\leq\tau$. We will show that $\BB_{\MM}^+$ is a base for $(\MM,\OO_u)$. Indeed, let $M\in\MM$, $\UU\in\OO_u$ and $M\in\UU$. Then $\exists U\in\BB_{\OO}$ such that $M\in U_{\MM}^+\subseteq\UU$. Thus $M\subseteq U$ and hence $\exists V\in\BB$ such that $M\subseteq V\subseteq U$. This implies that $M\in V_{\MM}^+\subseteq U_{\MM}^+$. Therefore, $M\in V_{\MM}^+\subseteq\UU$ and $V_{\MM}^+\in\BB_{\MM}^+$. Hence, $\BB_{\MM}^+$ is a base for $(\MM,\OO_u)$. Since $|\BB_{\MM}^+|\leq|\BB|\leq\tau$, we get that $w(\MM,\OO_u)\leq\tau$. Having in mind that $\OO=\OO_u\vee\OO_l$ (see Fact \ref{supremum}), we obtain that $w(\MM,\OO)\leq\tau$. \sqs

The space $(Comp(X),(\UP_{X})|_{Comp(X)})$ will be denoted by $\ZZ(X)$.

\begin{cor}\label{10'}{\rm (\cite{Mi})}
Let $(X,\TT)$ be a $T_1$-space and $w(X)\geq\aleph_0$. Then $w(\ZZ(X))=w(X)$.
\end{cor}

\doc Let $\OO$ be the restriction of the Vietoris topology $\UP_X$ on $Comp(X)$. Then $\ZZ(X)=(Comp(X),\OO)$. Using Proposition \ref{strV}, we get that $\OO$ is a strong Vietoris-type topology on $Comp(X)$ and $\TT=\TT_{\OO}$. Hence $\PP_{\OO}=\TT=\BB_{\OO}$. Then, clearly, $w(X,\PP_{\OO})=w(X)$. Let $\BB_0$ be a base for $(X,\TT)$ such that $|\BB_0|=w(X)$. Then $\BB=(\BB_0)^{\cup}$ satisfies the hypothesis of Proposition \ref{teglovi} for $\tau=|\BB|=w(X)$. Hence, by Proposition \ref{teglovi}, $w(\ZZ(X))\leq w(X)$. Since $X$ can be embedded in $\ZZ(X)$ (by Proposition \ref{2.7.20b}), we have that $w(X)\leq w(\ZZ(X))$. Therefore, $w(X)=w(\ZZ(X))$. \sqs

\begin{pro}\label{2.7.20e}
Let $(X,\TT)$ be a space, $\MM\subseteq\PP'(X)$, $Fin(X)\subseteq\MM$, $\OO$ be a Vietoris-type topology on $\MM$ and $\TT_{\OO}=\TT$. If $d(X)\geq\aleph_0$ then $d(X)\geq d(\MM,\OO)$.
\end{pro}

\doc Let $d(X)=\tau\geq\aleph_0$. Then there exists an $A\subset X$ such that $\overline{A}=X$ and $|A|=\tau$. Hence $|Fin(A)|=\tau$. We will prove that $Fin(A)$ is dense in $(\MM,\OO)$. Indeed, for every $\displaystyle U_{\MM}^+\cap\bigcap_{i=1}^n(V_i)_{\MM}^-=O\not=\emptyset$, where $U\in\BB_{\OO}$ and $V_i\in\PP_{\OO}$ for $i=1,\dots,n$, there exist $x_i\in A\cap V_i\cap U$. Hence $\{x_i\st i=1,\dots,n\}\in Fin(A)\cap O$. Thus $d(\MM,\OO)\leq\tau$. Therefore $d(\MM,\OO)\leq d(X)$. \sqs

\begin{pro}\label{T_0}
Let $X$ be a set, $\MM\subseteq\PP\ap(X)$ and $\OO$ be a topology on $\MM$. If $\MM\subseteq\{X\setminus A\st A\in\PP_{\OO}\}$ or $\MM\subseteq\BB_{\OO}$ then $(\MM,\OO)$ is a $T_0$-space.
\end{pro}

\doc Let $M,M_1\in\MM$ and $M\not=M_1$. Then $M\setminus M_1\not=\emptyset$ or $M_1\setminus M\not=\emptyset$. Let $M_1\setminus M\not=\emptyset$. If $\MM\subseteq\{X\setminus A\st A\in\PP_{\OO}\}$ then, setting $O=(X\setminus M)_{\MM}^-$, we get that $O\in\OO$, $M_1\in O$ and $M\not\in O$. If $\MM\subseteq\BB_{\OO}$ then, setting $O\ap=M_{\MM}^+$, we get that $O\ap\in\OO$, $M\in O\ap$ and $M_1\not\in O\ap$. We argue analogously if $M\setminus M_1\not=\emptyset$. \sqs

\begin{cor}\label{cT_0}{\rm (\cite{Mi})}
If $(X,\TT)$ is a topological space then $(CL(X),\UP_X)$ is a $T_0$-space.
\end{cor}

\doc Set $\OO=\UP_X$. Then $\TT\subseteq\PP_{\OO}$. Hence $\MM=CL(X)\subseteq\{X\setminus A\st A\in\PP_{\OO}\}$. Thus, by Proposition \ref{T_0}, $(\MM,\OO)$ is a $T_0$-space. \sqs

\begin{pro}\label{T_1}
Let $(X,\TT)$ be a topological space, $\MM\subseteq\PP'(X)$, $\PP$ be a subbase for $(X,\TT)$ and $\MM\subseteq\{X\setminus U\st U\in\PP\}$. Let $\MM$ be a natural family and $\OO$ be the topology on $\MM$ having as a subbase the family $\PP_{\MM}^-\cup\PP_{\MM}^+$. Then $\OO$ is a strong Vietoris-type topology on $\MM$, $\TT_{\OO}=\TT$ and $(\MM,\OO)$ is a $T_1$-space.
\end{pro}

\doc Since $\bigcup\PP=X$, we get easily that $\bigcup\PP_{\MM}^-=\MM$. Thus $\PP_{\MM}^-\cup\PP_{\MM}^+$ can serve as a subbase for a topology on $\MM$. Clearly, $\OO$ is a Vietoris-type topology on $\MM$, $\PP\subseteq\PP_{\OO}$ and $\PP^{\cap}\subseteq\BB_{\OO}$. We will show that $\PP_{\OO}\cup\BB_{\OO}\subseteq\TT$. Indeed, let $A\in\BB_{\OO}$. Let $x\in A$. Then $\{x\}\in A_{\MM}^+\in\OO$ and hence $\exists U\in\PP^{\cap}$ and $\exists V_1,\dots,V_n\in\PP$ such that $\displaystyle\{x\}\in U_{\MM}^+\cap\bigcap_{i=1}^n(V_i)_{\MM}^-=O\subseteq A_{\MM}^+$. Let $\displaystyle V=\bigcap_{i=1}^nV_i$. Then $x\in U\cap V\in\PP^{\cap}\subseteq\TT$. We will show that $U\cap V\subseteq A$. Indeed, let $y\in U\cap V$. Then $\{y\}\in O\subseteq A_{\MM}^+$ and thus $y\in A$. Therefore, $x\in U\cap V\subseteq A$ and $U\cap V\in\TT$. Hence, $A\in\TT$. Analogously, we get that if $A\in\PP_{\OO}$ then $A\in\TT$. Hence $\PP\subseteq\PP_{\OO}\cup\BB_{\OO}\subseteq\TT$. This implies that $\TT_{\OO}=\TT$.

Since we also have that $\PP\subseteq\PP_{\OO}\subseteq\TT$ and $\PP\subseteq\BB_{\OO}\subseteq\TT$, we get that $\TT_{-\OO}=\TT=\TT_{+\OO}$. Thus $\OO$ is a strong Vietoris-type topology on $\MM$.

Let $M\in\MM$. We will prove that $\overline{\{M\}}^{\OO}=\{M\}$. Indeed, let $M_1\in\MM$ and $M_1\not=M$. Then $M_1\setminus M\not=\emptyset$ or $M\setminus M_1\not=\emptyset$. Let $M_1\setminus M\not=\emptyset$. Then $U=X\setminus M\in\PP$, $M_1\in U_{\MM}^-$ and $M\not\in U_{\MM}^-$. Hence $M_1\not\in\overline{\{M\}}^{\OO}$. Let now $M\setminus M_1\not=\emptyset$. Let $x\in M\setminus M_1$ and $U=X\setminus\{x\}$. Then $U\in\PP$, $M_1\in U_{\MM}^+$ and $M\not\in U_{\MM}^+$. Hence $M_1\not\in\overline{\{M\}}^{\OO}$. Therefore $\overline{\{M\}}^{\OO}=\{M\}$. So, $(\MM,\OO)$ is a $T_1$-space. \sqs

\begin{rem}\label{3.14'}
\rm
The above proof shows that the condition ''$\MM\subseteq\{X\setminus U\st U\in\PP\}$'' in Proposition \ref{T_1} can be replaced by the following one:

\smallskip

\noindent (*)\ \ \ if $M,M\ap\in\MM$ and $M\setminus M\ap\not=\emptyset$ then $\exists U,V\in\PP$ such that $M\ap\subseteq U$, $M\not\subseteq U$, $M\cap V\not=\emptyset$ and $M\ap\cap V=\emptyset$.

\smallskip

\noindent Note that condition (*) is fulfilled if the following condition holds:

\smallskip

\noindent (**)\ \ \ if $M\in\MM$ and $x\not\in M$ then $\exists U,V\in\PP$ such that $M\subseteq U\subseteq X\setminus\{x\}$ and $x\in V\subseteq X\setminus M$.
\end{rem}

\begin{cor}\label{T_1cor}{\rm (\cite{Mi})}
Let $(X,\TT)$ be a $T_1$-space. Then $(CL(X),\UP_X)$ is a $T_1$-space.
\end{cor}

\doc Set $\PP=\TT$ and $\MM=CL(X)$ in Proposition \ref{T_1}. \sqs

\begin{defi}\label{P-reg}
\rm
Let $(X,\TT)$ be a topological space and $\PP$ be a subbase for $(X,\TT)$. $(X,\TT)$ is said to be $\PP$-{\em regular}, if for every $x\in X$ and for every $U\in\PP$ such that $x\in U$, there exist $V,W\in\PP$ with $x\in V\subset X\setminus W\subset U$.
\end{defi}

Clearly, a topology space $(X,\TT)$ is $\TT$-regular iff it is regular. Also, every $\PP$-regular space is regular.

\begin{exa}\label{notPreg}
\rm
Let $(\RRRR,\TT)$ be the real line with its natural topology and $\PP=\{(\alpha,\beta)\setminus F\st\alpha,\beta\in\overline{\RRRR},\ \alpha<\beta,\ F\subset\RRRR,\ |F|<\aleph_0\}$ or $\PP=\{(\alpha,\beta)\st\alpha,\beta\in\overline{\RRRR},\ \alpha<\beta\}$. Then $\PP$ is a base for $\TT$, $\PP^{\cap}=\PP$, $(\RRRR,\TT)$ is regular but it is not $\PP$-regular.
\end{exa}

\begin{pro}\label{T_2}
Let $(X,\TT)$ be a topological space, $\MM\subseteq\PP'(X)$, $\PP$ be a subbase for $(X,\TT)$ and $\MM\subseteq\{X\setminus U\st U\in\PP\}$. Let $\MM$ be a natural family, $X$ be $\PP$-regular and $\OO$ be the topology on $\MM$ having as a subbase the family $\PP_{\MM}^-\cup\PP_{\MM}^+$. Then $\OO$ is a strong Vietoris-type topology on $\MM$, $\TT_{\OO}=\TT$ and $(\MM,\OO)$ is a $T_2$-space.
\end{pro}

\doc From Proposition \ref{T_1}, we get that $\OO$ is a strong Vietoris-type topology on $\MM$ and $\TT_{\OO}=\TT$. Let us prove that $(\MM,\OO)$ is a $T_2$-space. Let $M_1\not=M_2$. Then $M_1\setminus M_2\not=\emptyset$ or $M_2\setminus M_1\not=\emptyset$. Let, for example, $M_1\setminus M_2\not=\emptyset$. Let us fix $x\in M_1\setminus M_2$. Then $U=X\setminus M_2\in\PP$ and $x\in U$. Hence, there exist $V,W\in\PP$, such that $x\in V\subset X\setminus W\subset U$. Then $V\cap W=\emptyset$ and $W\supset X\setminus U=M_2$. Therefore $M_1\in V_{\MM}^-$, $M_2\in W_{\MM}^+$ and $V_{\MM}^-\cap W_{\MM}^+=\emptyset$. If $M_2\setminus M_1\not=\emptyset$ then we argue analogously. \sqs

\begin{pro}\label{proPreg}
Let $(X,\TT)$ be a topological space, $\MM\subseteq\PP\ap(X)$, $\PP$ be a base for $(X,\TT)$, $\MM=\{X\setminus U\st U\in\PP\}$ and $\PP=\PP^{\cap}$. Let $\MM$ be a natural family, $\OO$ be the topology on $\MM$ having as a subbase the family $\PP_{\MM}^-\cup\PP_{\MM}^+$ and $(\MM,\OO)$ be a $T_2$-space. Then $X$ is $\PP$-regular.
\end{pro}

\doc Note that $Fin(X)\subseteq\MM$ and $\MM=\MM^{\cup}$. Let now $U\in\PP$ and $x\in U$. Then $\{x\}\in\MM$, $X\setminus U\in\MM$, $\{x\}\cup(X\setminus U)\in\MM$ and $\{x\}\cup(X\setminus U)\not=X\setminus U$. Since $(\MM,\OO)$ is a $T_2$-space, there exist $V,W,V_1,\dots,V_n,W_1,\dots,W_k\in\PP$ such that $\displaystyle\{x\}\cup(X\setminus U)\in V_{\MM}^+\cap\bigcap_{i=1}^n(V_i)_{\MM}^-=O$, $\displaystyle X\setminus U\in W_{\MM}^+\cap\bigcap_{j=1}^k(W_j)_{\MM}^-=O\ap$ and $O\cap O\ap=\emptyset$. Then $\{x\}\cup(X\setminus U)\subseteq V$, $(\{x\}\cup(X\setminus U))\cap V_i\not=\emptyset$ for $i=1,\dots,n$, $X\setminus U\subseteq W$ and $(X\setminus U)\cap W_j\not=\emptyset$ for $j=1,\dots,k$. Fix $x_j\in(X\setminus U)\cap W_j$ for every $j=1,\dots,k$. Set $F=\{x_j\st j=1,\dots,k\}$. If $(X\setminus U)\cap V_i\not=\emptyset$ for $i=1,\dots,n$ then for each $i\in\{1,\dots,n\}$ there exists $y_i\in(X\setminus U)\cap V_i$. Then $F\cup\{y_i\st i=1,\dots,n\}\in O\cap O\ap$, a contradiction. Hence, there exists $i_0\in\{1,\dots,n\}$ such that $(X\setminus U)\cap V_{i_0}=\emptyset$. Then $x\in V_{i_0}$. Let $G=\{i\in\{1,\dots,n\}\st(X\setminus U)\cap V_i=\emptyset\}$. Then $G\not=\emptyset$ and $\displaystyle x\in\bigcap_{i\in G}V_i=V\ap$. Let $H=\{1,\dots,n\}\setminus G$. If $H\not=\emptyset$ then, for every $i\in H$, fix $z_i\in(X\setminus U)\cap V_i$ and set $F\ap=\{z_i\st i\in H\}$. If $H=\emptyset$ then set $F\ap=\emptyset$. We will show that $V\cap V\ap\cap W=\emptyset$. Indeed, suppose that there exists $y\in V\cap V\ap\cap W$ and set $F''=F\cup F\ap\cup\{y\}$. Then $F''\in O\cap O\ap$, a contradiction. Hence $V''=V\cap V\ap\subseteq X\setminus W$. Therefore, we obtain that $V''\in\PP$ and $x\in V''\subseteq X\setminus W\subseteq U$. Thus, $X$ is $\PP$-regular. \sqs

\begin{rem}\label{3.16''}
\rm
The above proof shows that the conditions $``\MM=\{X\setminus U\st U\in\PP\}$" and $``\MM$ is a natural family" in Proposition \ref{proPreg} can be replaced by $``\MM\supseteq\{X\setminus U\st U\in\PP\}$, $Fin X\subseteq\MM$ and if $x\in U\in\PP$ then $\{x\}\cup(X\setminus U)\in\MM$".
\end{rem}

\begin{cor}\label{T_3}{\rm (\cite{Mi})}
Let $(X,\TT)$ be a $T_1$-space. Then $(X,\TT)$ is a $T_3$-space iff $(CL(X),\UP_X)$ is a $T_2$-space.
\end{cor}

\doc Set $\PP=\TT$ and $\MM=CL(X)$ in Propositions \ref{T_2} and \ref{proPreg}. \sqs

\begin{pro}\label{mycom}
Let $(X,\TT)$ be a compact $T_1$-space, $\OO$ be a Vietoris-type topology on $CL(X,\TT)$ and $\TT_{\OO}=\TT$. Then $(CL(X,\TT),\OO)$ is a compact space.
\end{pro}

\doc Clearly, the identity map $i:(CL(X,\TT),\UP_X)\longrightarrow(CL(X,\TT),\OO)$ is continuous. By a theorem of E. Michael \cite[Theorem 4.2]{Mi}, the space $(CL(X,\TT),\UP_X)$ is compact. Hence  $(CL(X,\TT),\OO)$ is compact. \sqs

\section{Subspaces and hyperspaces}

In this section, we will regard the problem of continuity or inverse continuity of the maps of the form $i_{A,X}$ (see Proposition \ref{i_A,Xnepr} below for the definition of the maps $i_{A,X}$) for the hyperspaces with a srong Vietoris-type topology. This problem was regarded by H.-J.Schmidt \cite{Sch} for the lower Vietoris topology, by G. Dimov \cite{D1,D2} for the (upper) Vietoris topology and by Barov- Dimov- Nedev \cite{BDN1,BDN2} for the upper Vietoris topology.

We will need the following result from \cite{ED}:

\begin{pro}\label{i_A}{\rm (\cite{ED})}
Let $(X,\TT)$ be a space, $\PP$ be a subbase for $\TT$, $X\in\PP$ and $\PP^-$ be a subbase for a topology $\OO$ on $CL(X)$. Let $A$ be a subspace of $X$. Set $\PP_A=\{U\cap A\st U\in\PP\}$. Let $\OO_-^A$ be the topology on $CL(A)$ having $(\PP_A)_{CL(A)}^-$ as a subbase. Then $i_{A,X,-}:(CL(A),\OO_-^A)\longrightarrow(CL(X),\OO)$, where $i_{A,X,-}(F)=\overline{F}^X$, is a homeomorphic embedding.
\end{pro}

The next assertion is trivial.

\begin{pro}\label{2.6^0}
Let $A$ and $X$ be sets, and $f:A\longrightarrow X$ be a function. Let, for $i=1,2$, $\TT_i$ (resp., $\OO_i$) be a topology on $A$ (resp., $X$). Let the maps $f:(A,\TT_1)\longrightarrow(X,\OO_1)$ and $f:(A,\TT_2)\longrightarrow(X,\OO_2)$ be continuous. Then $f:(A,\TT_1\vee\TT_2)\longrightarrow(X,\OO_1\vee\OO_2)$ is a continuous map.
\end{pro}

\begin{pro}\label{i_A,Xnepr}
Let $(X,\TT)$ be a $T_1$-space, $\MM=CL(X,\TT)$, $\OO$ be a strong Vietoris-type topology on $\MM$ and $\TT_{\OO}=\TT$. Let $A$ be a subspace of $X$, $\MM_A=CL(A)$, $\OO_-^A$ be the topology on $\MM_A$ having as a subbase the family $(\PP_{\OO}^A)_{CL(A)}^-$, where $\PP_{\OO}^A=\{U\cap A\st U\in\PP_{\OO}\}$, and let $\OO_+^A$ be the topology on $\MM_A$ having as a base the family $(\BB_{\OO}^A)_{CL(A)}^+$, where $\BB_{\OO}^A=\{U\cap A\st U\in\BB_{\OO}\}$. Then:

\smallskip

\noindent (a) $\OO^A=\OO_-^A\vee\OO_+^A$ is a strong Vietoris-type topology on $\MM_A$, and

\smallskip

\noindent (b) the map $i_{A,X}:(\MM_A,\OO^A)\longrightarrow(\MM,\OO)$, where $i_{A,X}(F)=\overline{F}^X$ for every $F\in\MM_A$, is continuous (resp., inversely continuous) if and only if the map $i_{A,X,+}:(\MM_A,\OO_+^A)\longrightarrow(\MM,\OO_u)$ is continuous (resp., inversely continuous) (here $i_{A,X,+}(F)=\overline{F}^X$ for every $F\in CL(A)$).
\end{pro}

\doc (a) Since $\TT_{+\OO}\equiv\TT_{-\OO}\equiv\TT_{\OO}\equiv\TT$, we get that $\PP_{\OO}$ is a subbase for $\TT$ and $\BB_{\OO}$ is a base for $\TT$. Also, $X\in\PP_{\OO}\cap\BB_{\OO}$. Then we easily obtain that $\TT_{+\OO^A}=\TT_{-\OO^A}=\TT_A$, where $\TT_A=\{U\cap A\st U\in\TT\}$. Hence $\OO^A$ is a strong Vietoris-type topology on $\MM_A$ and $\TT_{\OO^A}=\TT_A$.

\smallskip

\noindent (b) Clearly, the map $i_{A,X}$ is an injection.

\smallskip

\noindent $(\Leftarrow)$ It follows from Propositions \ref{i_A} and \ref{2.6^0}.

\smallskip

\noindent $(\Rightarrow)$ Let $i_{A,X}:(\MM_A,\OO^A)\longrightarrow(\MM,\OO)$ be continuous. We will prove that $i_{A,X,+}:(\MM_A,\OO_+^A)\longrightarrow(\MM,\OO_u)$ is continuous. Let $U\in\BB_{\OO}$, $F\in\MM_A$ and $\overline{F}^X\subset U$. Then there exist $V\in\BB_{\OO}^A$ and $V_1,\dots,V_n\in\PP_{\OO}^A$ such that $\displaystyle F\in V_{\MM_A}^+\cap\bigcap_{i=1}^n(V_i)_{\MM_A}^-=O$ and $i_{A,X}(O)\subseteq U_{\MM}^+$. Thus $F\cap V_i\not=\emptyset$, $\forall i=1,\dots,n$ and $F\subseteq V$. Let $a_i\in F\cap V_i$, for $i=1,\dots,n$. We will show that $i_{A,X,+}(V_{\MM_A}^+)\subset U_{\MM}^+$. Indeed, let $G\in\MM_A$ and $G\subset V$. Set $G\ap=G\cup\{a_i\st i=1,\dots,n\}$. Since $X$ is a $T_1$-space, we get that $G\ap\in\MM_A$. Obviously, $G\ap\in O$. Hence $\overline{G\ap}^X\subset U$. Thus  $\overline{G}^X\subset U$. So, $i_{A,X,+}(V_{\MM_A}^+)\subseteq U_{\MM}^+$. Therefore, the map $i_{A,X,+}$ is continuous.

Let $i_{A,X}$ be inversely continuous. We will show that  the map $i_{A,X,+}$ is inversely continuous. Let $U\in\BB_{\OO}^A$, $F\in\MM_A$ and $F\subset U$. Then there exist $V\in\BB_{\OO}$ and $V_1,\dots,V_n\in\PP_{\OO}$ such that $\displaystyle\overline{F}^X\in V^+\cap\bigcap_{i=1}^nV_i^-=O$ and $(i_{A,X})^{-1}(O)\subseteq U_{\MM_A}^+$. Thus $\overline{F}^X\subset V$ and $\overline{F}^X\cap V_i\not=\emptyset$, $\forall i=1,\dots,n$. Obviously, $\forall i=1,\dots,n$, there exists $a_i\in F\cap V_i$. We will show that $(i_{A,X,+})^{-1}(V^+)\subseteq U_{\MM_A}^+$. Indeed, let $G\in\MM_A$ and $\overline{G}^X\subset V$. Then $G\ap=G\cup\{a_i\st i=1,\dots,n\}\in\MM_A$ and $\overline{G\ap}^X\in O$. Hence $G\ap\subset U$. Thus $G\subset U$. So, $(i_{A,X,+})^{-1}(V^+)\subseteq U_{\MM_A}^+$. Therefore, the map $i_{A,X,+}$ is inversely continuous. \sqs

\baselineskip = 0.75\normalbaselineskip

\end{document}